# A NEW PROOF OF GOODSTEIN'S THEOREM

JUAN A. PEREZ

ABSTRACT. Goodstein sequences are numerical sequences in which a natural number $m$, expressed as the complete normal form to a given base $a$, is modified by increasing the value of the base $a$ by one unit and subtracting one unit from the resulting expression. As initially defined, the first term of the Goodstein sequence is the complete normal form of $m$ to base 2. Goodstein's Theorem states that, for all $m \in \mathbb{N}$, the Goodstein sequence eventually terminates at zero. Goodstein's Theorem was originally proved using the well-ordered properties of transfinite ordinals. The theorem was also shown to be unprovable-in-PA (Peano Arithmetic) using transfinite induction and Gödel's Second Incompleteness Theorem.

This article describes a proof of Goodstein's Theorem in first-order arithmetic that contradicts the theorem's unprovability-in-PA. The proof uses mathematical induction and is applied (via the super-exponential function) to a generalized version of the Goodstein sequences. Such a proof demonstrates the inconsistency of classical set theory, more precisely the combination of the Zermelo-Fraenkel axioms and the axiom of choice (ZFC).

## 1. INTRODUCTION

The concept of mathematical incompleteness was first brought to universal attention by Gödel's Theorems [3,6,10]. The initial search for new examples of incompleteness produced statements that were basically Gödelian in character, i.e. 'true-but-unprovable' sentences of a logical nature, constructed using coding devices [10]. Therefore, it was important to discover other more natural sentences that, as well as being of more intrinsic mathematical interest, could be shown to be both true and unprovable in first-order arithmetic (commonly Peano Arithmetic, PA). The verification of the truth and unprovability of such examples has to depend on mathematical tools that, while unavailable to basic arithmetics, are also non-Gödelian in character. The published examples achieve this goal by making use of transfinite numbers [10].

A case of particular interest and significance is Goodstein's Theorem [5,9,10]. This theorem has a purely number-theoretic character and deals with basic arithmetic sequences. It was originally proved (by Goodstein, in 1944), not by means of first-order arithmetic, but by making use of the well-ordered properties of transfinite ordinals [4,5]. It was later shown (by Paris and Kirby, in 1982) to be unprovable-in-PA, with a proof that needs both transfinite induction and Gödel's

Second Incompleteness Theorem [8,10]. Paris and Kirby's argument requires a previous result, i.e. the proof by Gentzen of the consistency of PA, also achieved by transfinite induction [1,2]. Combining all these elements, Paris and Kirby showed that, if PA were capable of proving Goodstein's Theorem, PA would also be able to prove its own consistency, thus violating Gödel's Second Theorem [8]. Therefore, the claim is made that PA cannot prove Goodstein's Theorem.

The results by Paris and Kirby establish a close link between the consistency of classical set theory, i.e. the combination of the Zermelo-Fraenkel axioms and the axiom of choice (ZFC) [4,7], and the inability of first-order arithmetic to prove Goodstein's Theorem. The implication is that the existence of an elementary proof of Goodstein's Theorem (if such a proof could be constructed) would both contradict Paris and Kirby's conclusion [8] and strongly challenge the integrity of ZFC [4,7]. This article describes a proof-in-PA of Goodstein's Theorem that uses mathematical induction and is applied (via the super-exponential function) to a generalized version of the theorem.

## 2. Goodstein sequences

In its original formulation, Goodstein's Theorem refers to numerical sequences that start from a natural number $m$, expressed as pure representation in base $a$ (also known as the complete normal form of $m$ to base $a$), where no number greater than $a$ can appear [5,9,10]. For example, if $m = 311$ and $a = 2$, the complete normal form of 311 to base 2 would be:

$$(2.1) \qquad 311 = 2^8 + 2^5 + 2^4 + 2^2 + 2^1 + 2^0 = 2^{2^{2+1}} + 2^{2^2+1} + 2^{2^2} + 2^2 + 2^1 + 2^0.$$

Similarly, if $a = 3$, the complete normal form of 311 to base 3 would be:

$$(2.2) \qquad 311 = 3^5 + 2 \cdot 3^3 + 3^2 + 3^1 + 2 \cdot 3^0 = 3^{3+2} + 2 \cdot 3^3 + 3^2 + 3^1 + 2 \cdot 3^0.$$

As the examples (2.1) and (2.2) illustrate, a complete normal form to a given base generates a summation of power terms, where each power can be multiplied by a coefficient always smaller than the base. It is immediately apparent that these coefficients are smaller than the base, since coefficients greater than the base are assimilated into powers with a greater exponent (for example, $4 \cdot 3^3 = 3 \cdot 3^3 + 1 \cdot 3^3 = 3^4 + 3^3$). The exponents of the powers are also represented as summation of power terms (e.g. $2^8 = 2^{2^{2+1}}$), so the process continues with as many levels of exponentiation as required. Each number $m \in \mathbb{N}$ has a unique pure representation in base $a$, whatever the value of $a \in \mathbb{N}$.

To generate the original Goodstein sequence, once the starting number $m$ is expressed as a pure representation in base $a$ (in the original version of the theorem, the initial value of $a$ is always 2), the next step is to produce a new number by replacing the base $a$ with a new base $a+1$, evaluating the result and then subtracting one. The resulting number (expressed as complete normal form to base $a+1$) then becomes the template for the next step, where a similar operation is undertaken, with the process continuing in the same way. This defines a set of primitive recursive functions denoted $F_n(r)$, where the number $r$ is first represented in base $n$ and then transformed as described [9,10]. These functions provide the basis of the Goodstein sequence $g(m,n)_{n\geq 2}$, which can be recursively defined as

$$g(m,2) = m \tag{2.3}$$
$$g(m,n+1) = F_n\big(g(m,n)\big) .$$

According to (2.3), the Goodstein sequence starting from $m$ is

$$m, F_2(m), F_3(F_2(m)), F_4(F_3(F_2(m))), \ldots$$

By starting the sequence with $g(m,2) = m$, the index $n$ of each term $g(m,n)$ coincides with the base of the resulting representation.

For numbers $m \geq 4$, the values of the functions $F_n(r)$ become extremely large very quickly, so it seems that the Goodstein sequences should grow indefinitely [9,10]. Yet all Goodstein sequences terminate at zero, as proved by Goodstein's Theorem [5,10]:

**Theorem** *For all $m \in \mathbb{N}$, the Goodstein sequence eventually terminates at zero:*
$$\forall m \in \mathbb{N},\ \exists n \in \mathbb{N}\ (g(m,n)=0).$$

The original formulation of Goodstein's Theorem accounts only for the simplest type of Goodstein sequences, i.e. those where the starting base is $a = 2$ and the functions $F_n(r)$ transform the corresponding base by increasing its value only one unit. However, a more general version of Goodstein sequences can be defined whereby the starting base can take any arbitrary value $a = a_0$, and the step increase in the value of the base applied by each function $F_n(r)$ is equally arbitrary (although finite). The objective of this paper is to prove by ordinary arithmetic that such generalized Goodstein sequences terminate at zero. Then, the conventional Goodstein's Theorem described becomes an example of the more general theorem.

## 3. Prerequisites of the theorem

Prior to the theorem, there are a few elementary results that need to be stated. First, it needs to be emphasized that the terms of a Goodstein sequence, for any finite numbers of steps, are also finite in value, regardless of the degree of repeated exponentiation. Thus the following lemma and corollary.

**Lemma 3.1.** *For any two numbers $a, b \in \mathbb{N}$, the operation of exponentation yields a new number:*
$$\forall a,b \in \mathbb{N},\ \exists c \in \mathbb{N}\ (c = a^b).$$

*Proof.* Since *i)* $a^1 = a$, $\forall a \in \mathbb{N}$, and *ii)* $1^b = 1$, $\forall b \in \mathbb{N}$, assume that $a, b \geq 2$. Consider $\mathbb{N}$, the set of the natural numbers. Now consider, within the enumeration of the elements of $\mathbb{N}$, all the powers of $a$, which can be selected to construct the corresponding subset $N_a \subseteq \mathbb{N}$:

$$N_a = \{a^0, a^1, a^2, a^3, a^4, \cdots\} = \{a^k, k \in \mathbb{N}\}. \tag{3.1}$$

The critical requirement is to establish that the subset $N_a$ is an infinite set. If $N_a$ were a finite set, the implication would be that one of the elements of $N_a$ would be greater than all the others:

$$N_a \text{ finite} \Rightarrow \exists l \in \mathbb{N} \left(\forall a^k \in N_a,\ a^k \leq a^l\right).$$

The product of any two natural numbers is also a natural number; therefore, it can be stated that the product $a \cdot a^l = a^{l+1}$ is an element of $\mathbb{N}$. However, the definition of $N_a$ establishes that any power of $a$ is an element of $N_a$. Accordingly, $a^{l+1} \in N_a$. Since $a^l < a^{l+1}$, this contradicts that $a^l$ is the greatest element in the subset; hence, $N_a$ cannot be finite. (3.1) can then be modified to accommodate the infinite nature of $N_a$,

(3.2) $$N_a = \{a^0, a^1, a^2, a^3, a^4, \cdots\} = \{a^k, \forall k \in \mathbb{N}\}.$$

From (3.2), it can be deduced that $a^b \in N_a$. Since $N_a \subseteq \mathbb{N}$, the final conclusion is reached that $a^b \in \mathbb{N}$. □

When analysing Goodstein sequences, a critical factor to consider is that the representation of $m$ as complete normal form to base $a$ can require the extraction of powers in base $a$ of any exponents greater than $a$; the same can apply to exponents within the exponents. Therefore, the possible presence of multiple levels of exponentiation needs to be addressed in a systematic way. This requires the introduction of the super-exponential function [10], which can be recursively defined as

(3.3) $$\begin{aligned} a \Uparrow 0 &= a \\ a \Uparrow Sn &= a^{a \Uparrow n}. \end{aligned}$$

For example, as defined by (3.3), $a \Uparrow 3$ will be $a^{a^{a^a}}$, with a 'tower' of 3 levels of exponentiation. It is important to establish that, provided the number of levels of exponentiation is finite, the resulting number is also finite. The purpose of the very basic Lemma 3.1 is to facilitate the formulation of a corollary that fulfils that requirement.

**Corollary 3.2.** *For any two numbers $a, n \in \mathbb{N}$, the super-exponential function yields a new number:*

$$\forall a, n \in \mathbb{N},\ \exists b \in \mathbb{N}\ (b = a \Uparrow n).$$

*Proof.* By repeated application of Lemma 3.1. Consider first the single exponential $b_1 = a^a$. Lemma 3.1 implies that $b_1 \in \mathbb{N}$. Now consider a new exponential $b_2 = a^{b_1}$. Since $b_1 \in \mathbb{N}$, Lemma 3.1 implies that $b_2 \in \mathbb{N}$. Repeating the same pattern a total of $n$ times, a final exponential $b_n = a^{b_{n-1}} = a \Uparrow n$ will be obtained that, because of Lemma 3.1, is also a natural number $b_n \in \mathbb{N}$. □

As indicated, the significance of Corollary 3.2 lies in the fact that, however large the value of $a$, the super-exponential $a \Uparrow n$ is always a finite number, provided $n$ is itself finite. This is relevant to all Goodstein sequences because, whatever the increase in value of the base used for each term of the sequence, there is never an increase in the number of levels of exponentiation anywhere in the representations; on the contrary, the number of levels of exponentiation gradually decreases. Consequently, for any finite number of steps, the value of the corresponding term for all Goodstein sequences is always finite. This implies that, provided a satisfactory method could be devised to analyse any Goodstein sequence, such an analysis would remain in first-order arithmetic, as would any subsequent proof of Goodstein's Theorem.

Another basic but fundamental prerequisite for the proof construction is provided by the following lemma.

**Lemma 3.3.** *For any two numbers $a, b \in \mathbb{N}$, it is always possible to construct a super-exponential $a \Uparrow n$ greater or equal than $b$:*

$$\forall a, b \in \mathbb{N}, \exists n \in \mathbb{N} \, (b \leq a \Uparrow n) \, .$$

*Proof.* Consider $\mathbb{N}$, the set of the natural numbers. Consider a subset $A_k \subseteq \mathbb{N}$ formed by all the super-exponentials $a \Uparrow k$ never greater than $b$:

$$A_k = \{a \Uparrow k, \, a, k \in \mathbb{N}, \, a \Uparrow k \leq b\} \, .$$

Since $A_k$ is a finite set, one of the elements of $A_k$ would be greater than all the others:

$$A_k \text{ finite} \Rightarrow \exists l \in \mathbb{N} \, (\forall a \Uparrow k \in A_k, \, a \Uparrow k \leq a \Uparrow l) \, .$$

Construct a new exponential $a \Uparrow (l+1) = a^{a \Uparrow l} > a \Uparrow l$, which is still a natural number (Lemma 3.1). Since it is greater than $a \Uparrow l$, the new exponential cannot be a member of $A_k$; therefore, it will be $a \Uparrow (l+1) = a \Uparrow n \geq b$, as required. □

The importance of Lemma 3.3 is made clear by the introduction of two further lemmas, which play a fundamental role in support of the final proof. Although Lemma 3.4 is a very elementary and well-known result, a proof is included for the sake of completeness.

**Lemma 3.4.** *For any two numbers $a, b \in \mathbb{N}$, $a, b \geq 1$, the operation of exponentiation $a^b$ followed by subtraction of one unit yields:*

$$\forall a, b \in \mathbb{N}, \, a, b \geq 1, \quad a^b - 1 = (a-1) \sum_{i=0}^{b-1} a^i \, .$$

*Proof.* By inspection. Rewriting the summation, we have that:

(3.4) $\qquad w = a^b - 1 = (a-1) \cdot a^{b-1} + (a-1) \cdot a^{b-2} + \cdots + (a-1) \cdot a^1 + (a-1) \cdot a^0 \, .$

The next step is to add back the subtracted unit to both sides of (3.4):

$$w + 1 = a^b = (a-1) \cdot a^{b-1} + (a-1) \cdot a^{b-2} + \cdots + (a-1) \cdot a^1 + (a-1) \cdot a^0 + 1$$
$$w + 1 = a^b = (a-1) \cdot a^{b-1} + (a-1) \cdot a^{b-2} + \cdots + (a-1) \cdot a^1 + a$$
$$w + 1 = a^b = (a-1) \cdot a^{b-1} + (a-1) \cdot a^{b-2} + \cdots + (a-1) \cdot a^2 + a^2$$
$$\cdots \cdots$$
$$w + 1 = a^b = (a-1) \cdot a^{b-1} + (a-1) \cdot a^{b-2} + a^{b-2}$$
$$w + 1 = a^b = (a-1) \cdot a^{b-1} + a^{b-1} = a^b \qquad \square$$

**Lemma 3.5.** *Consider two numbers $x, y \in \mathbb{N}$, $x, y \geq 1$. The Goodstein sequences of $x$ and $y$ always satisfy the condition*

$$\forall x, y, n \in \mathbb{N} \, (x > y \Rightarrow g(x, n) \geq g(y, n)) \, ,$$

*where the equal sign only applies once both sequences have terminated at zero.*

*Proof.* Elementary. Inspection of equation (3.4) helps to confirm this lemma, by comparing the increase in value that the various powers experience with each transformation undertaken by the functions $F_n(r)$; it is obvious that the powers with larger exponents will increase in value more than the powers with smaller exponents. Since the complete normal forms of $x$ and $y$ to base 2 that initiate the sequences are such that $x > y$, the powers that constitute $g(x, 2)$ are, on balance, greater in value than those that constitute $g(y, 2)$. Hence, the previous consideration ensures that the terms of $g(x, n)$ are always greater in value than the corresponding terms of $g(y, n)$.

The above conclusion is not in contradiction with the possibility of the terms of the two sequences equalizing, $g(x, n) = g(y, n)$, which will happen once both sequences have terminated at zero. For a Goodstein sequence to start decreasing in value, the powers left in the complete normal form can only be of a base to power zero, whatever the value of the coefficient that multiplies that base. The previous conclusion implies that $g(x, n)$ always takes longer than $g(y, n)$ to reach such a point, so $g(x, n) < g(y, n)$ can never occur. □

When considering Lemma 3.4, a key observation derived from equation (3.4) is that all the terms in the expression carry the same coefficient $(a-1)$, which is determined by the base $a$ but totally independent of the exponent $b$. This is of great significance for the construction of Goodstein sequences since, in all representations, the coefficients (already smaller than the base) remain largely unmodified by the transformations undertaken by the functions $F_n(r)$.[1] Furthermore, when $b \leq a$, the independence of the coefficients on subsequent increases in the value of the base also extends to the exponents, as can be deduced from equation (3.4).

The independence of the coefficients on the increases in base value, extended to the exponents when they no longer exceed the value of the base, and combined with the fact that the initial levels of exponentiation cannot increase within the same Goodstein sequence, are all key factors contributing to the eventual termination at zero of all Goodstein sequences.

Lemma 3.5 implies that, if the Goodstein sequence of a given number $m$ terminates at zero, the same is true of the Goodstein sequences of all other numbers smaller than $m$. This observation is a simple but important aid to the construction of a proof for Goodstein's Theorem; the proof can be simplified by the combination of Lemma 3.3 and Lemma 3.5 since the implication is that, in order to prove the whole theorem, it is sufficient to construct a proof for a reduced set of Goodstein sequences $g(m, n)$, as defined by (2.3), where $m = a \Uparrow k$, $\forall k \in \mathbb{N}$.

The final factor that confirms Goodstein's Theorem is that the termination at zero of all Goodstein sequences is independent of the actual increase in the value of the base, regardless of the size of such increases. As the construction of the proof will show, this independence plays a crucial role in the termination of all exponents, regardless of the number of levels of exponentiation.

---

[1] It is important to stress that the coefficients that multiply the base of the various powers that constitute each term of the Goodstein sequence are unaffected by subsequent transformations, with the only and obvious exception of the coefficient that multiplies the base to power zero.

## 4. Generalized Goodstein sequences

To accommodate variable increases in the value of the base, a generalized Goodstein sequence can be introduced, whereby the base of the initial complete normal form and the step increases in base values are both universal. On this basis, the sequence can be defined recursively as:

$$\begin{aligned} g_a(m, a_0) &= m\,,\ a_0 \geq 2 \\ g_a(m, a_{n+1}) &= F_{a_n, d_{n+1}}(g_a(m, a_n)) \\ a_{n+1} &= a_n + d_{n+1}\,,\ d_{n+1} \geq 1\,, \end{aligned} \tag{4.1}$$

where the two subindexes of the transformation function $F_{a_n, d_{n+1}}$ represent the initial base ($a_n$) for that step in the sequence, and the increase in the value of the base ($d_{n+1}$), respectively. Therefore, the generalized Goodstein sequence starting from $m$ is

$$m\,,\ F_{a_0, d_1}(m)\,,\ F_{a_1, d_2}(F_{a_0, d_1}(m))\,,\ F_{a_2, d_3}(F_{a_1, d_2}(F_{a_0, d_1}(m)))\,,\ \ldots$$

with

$$a_1 = a_0 + d_1\,,\ a_2 = a_1 + d_2\,,\ a_3 = a_2 + d_3\,,\ \ldots$$

Note that, when $a_0 = 2$ and $d_n = 1$, $\forall n \in \mathbb{N}$, the definition (4.1) readily simplifies to (2.3), the definition of the classical Goodstein sequence [5,9,10].

The property of conventional Goodstein sequences described by Lemma 3.5 will also hold for generalized Goodstein sequences.

**Lemma 4.1.** *Consider two numbers $x, y \in \mathbb{N}$, $x, y \geq 1$. The generalized Goodstein sequences of $x$ and $y$ always satisfy the condition*

$$\forall x, y, a_n \in \mathbb{N}\ (x > y \Rightarrow g_a(x, a_n) \geq g_a(y, a_n))\,,$$

*provided the initial base, $a_0$, and the step increases in base value used in both sequences are identical.*

*Proof.* See Lemma 3.5. The values of the initial base, $a_0$, and of the step increases in base value that generalized Goodstein sequences use are arbitrary. Therefore, the condition of equality required by the lemma can be easily met. □

## 5. Goodstein's theorem(s)

The resulting theorem can now be formulated:

**Theorem 5.1 (Generalized Goodstein's Theorem).** *For all $m \in \mathbb{N}$, the generalized Goodstein sequence eventually terminates at zero:*

$$\forall m \in \mathbb{N},\ \exists a_n \in \mathbb{N}\,(g_a(m, a_n) = 0).$$

*Proof.* This proof requires both inspection and induction. The combined use of Lemma 3.3 and Lemma 4.1 allows simplification of the proof to a reduced set of numbers, namely those resulting from the super-exponential function (3.3), i.e. $m = a \Uparrow k$, $\forall k \in \mathbb{N}$. One of the advantages of working with these numbers is that the presence of multiple levels of exponentiation is taken into account. Furthermore, it is possible to establish a clear induction process:

*Step i*) Prove that the sequence $g_a(m, a_n)$ terminates at zero when $m = a \Uparrow 1 = a^a$.

*Step ii*) Prove that, if the sequence $g_a(m, a_n)$ terminates at zero when $m = a \Uparrow k_i$, then it also terminates at zero when $m = a \Uparrow (k_i + 1)$.

*Step iii*) The verification of steps *i*) and *ii*) implies that the sequence $g_a(m, a_n)$ terminates at zero when $m = a \Uparrow k$, $\forall k \in \mathbb{N}$. This is true $\forall a \in \mathbb{N}$.

It is an important observation that, because of Lemma 3.4, the application of $F_{a_0, d_1}(m)$ to any of the super-exponentials $a \Uparrow k$ automatically reduces the levels of exponentiation by one unit. This highlights what is the recurrent theme in any Goodstein sequence, i.e. the progressive and relentless reduction in value of coefficients, exponents and levels of exponentiation in the complete normal forms that make up the terms of the sequence. Although somewhat tortuous, the systematic examination of this loss of value provides proofs for *i*) and *ii*).

*Step i*) $m = a \Uparrow 1 = a^a$. To construct the generalized Goodstein sequence, follow (4.1) and apply Lemma 3.4 repeatedly, starting from $m = a_0 \Uparrow 1 = a_0^{a_0}$:

$$g_a(m, a_1) = a_1^{a_1} - 1 = (a_1 - 1) \cdot a_1^{a_1 - 1} + (a_1 - 1) \cdot a_1^{a_1 - 2} + \cdots + (a_1 - 1) \cdot a_1^1 + (a_1 - 1) \cdot a_1^0,$$

where $a_1 = a_0 + d_1$. The next term will be

$$g_a(m, a_2) = (a_1 - 1) \cdot a_2^{a_1 - 1} + (a_1 - 1) \cdot a_2^{a_1 - 2} + \cdots + (a_1 - 1) \cdot a_2^2 + (a_1 - 1) \cdot a_2^1 + (a_1 - 1) \cdot a_2^0 - 1,$$

where $a_2 = a_1 + d_2$. Here, $(a_1 - 1) \cdot a_2^0 - 1 = (a_1 - 2) \cdot a_2^0$. Thus, the coefficients $(a_1 - 1)$ are unmodified and will remain unaltered in subsequents terms of the sequence, apart from the coefficient that multiplies the smallest power in the summation. Accordingly,

$$g_a(m, a_3) = (a_1 - 1) \cdot a_3^{a_1 - 1} + (a_1 - 1) \cdot a_3^{a_1 - 2} + \cdots + (a_1 - 1) \cdot a_3^2 + (a_1 - 1) \cdot a_3^1 + (a_1 - 2) \cdot a_3^0 - 1,$$

where $a_3 = a_2 + d_3$. The reduction in value of the last coefficient will continue until a term is reached such that

$$g_a(m, a_\alpha) = (a_1 - 1) \cdot a_\alpha^{a_1 - 1} + (a_1 - 1) \cdot a_\alpha^{a_1 - 2} + \cdots + (a_1 - 1) \cdot a_\alpha^2 + (a_1 - 1) \cdot a_\alpha^1 - 1,$$

where $a_\alpha = a_{\alpha - 1} + d_\alpha$. Here, $(a_1 - 1) \cdot a_\alpha^1 - 1 = (a_1 - 2) \cdot a_\alpha^1 + (a_\alpha - 1) \cdot a_\alpha^0$. Therefore,

$$g_a(m, a_{\alpha + 1}) = (a_1 - 1) \cdot a_{\alpha + 1}^{a_1 - 1} + (a_1 - 1) \cdot a_{\alpha + 1}^{a_1 - 2} + \cdots + (a_1 - 2) \cdot a_{\alpha + 1}^1 + (a_\alpha - 1) \cdot a_{\alpha + 1}^0 - 1,$$

where $a_{\alpha + 1} = a_\alpha + d_{\alpha + 1}$. Here, $(a_\alpha - 1) \cdot a_{\alpha + 1}^0 - 1 = (a_\alpha - 2) \cdot a_{\alpha + 1}^0$. As before, the decrease in value of the last coefficient continues unabated, until a term is reached such that

$$g_a(m, a_\beta) = (a_1 - 1) \cdot a_\beta^{a_1 - 1} + (a_1 - 1) \cdot a_\beta^{a_1 - 2} + \cdots + (a_1 - 1) \cdot a_\beta^2 + (a_1 - 2) \cdot a_\beta^1 - 1,$$

where $a_\beta = a_{\beta - 1} + d_\beta$. Here, $(a_1 - 2) \cdot a_\beta^1 - 1 = (a_1 - 3) \cdot a_\beta^1 + (a_\beta - 1) \cdot a_\beta^0$. Therefore,

$$g_a(m, a_{\beta + 1}) = (a_1 - 1) \cdot a_{\beta + 1}^{a_1 - 1} + (a_1 - 1) \cdot a_{\beta + 1}^{a_1 - 2} + \cdots + (a_1 - 3) \cdot a_{\beta + 1}^1 + (a_\beta - 1) \cdot a_{\beta + 1}^0 - 1,$$

where $a_{\beta + 1} = a_\beta + d_{\beta + 1}$. Similarly, $(a_\beta - 1) \cdot a_{\beta + 1}^0 - 1 = (a_\beta - 2) \cdot a_{\beta + 1}^0$. The pattern will repeat itself until the coefficient in $(a_1 - 3) \cdot a_{\beta + 1}^1$ is eliminated. The coefficients that multiply the last power in the term increase in value with each iteration (the latest one was $a_\beta - 2$ in $(a_\beta - 2) \cdot a_{\beta+1}^0$), so the process of elimination will require many terms. Nevertheless, such coefficients remain finite and unaffected by later changes in base. Consequently, the sequence will continue until a term is reached such that

$$g_a(m, a_\chi) = (a_1 - 1) \cdot a_\chi^{a_1 - 1} + (a_1 - 1) \cdot a_\chi^{a_1 - 2} + \cdots + (a_1 - 1) \cdot a_\chi^3 + (a_1 - 1) \cdot a_\chi^2 - 1,$$

where $a_\chi = a_{\chi-1} + d_\chi$. Since $(a_1 - 1) \cdot a_\chi^2 - 1 = (a_1 - 2) \cdot a_\chi^2 + a_\chi^2 - 1$, Lemma 3.4 means that $a_\chi^2 - 1 = (a_\chi - 1) \cdot a_\chi^1 + (a_\chi - 1) \cdot a_\chi^0$. This implies that (starting with base $a_\chi$ rather than base $a_1$) the tail of the last term, i.e. $a_\chi^2 - 1$, repeats the same pattern previously encountered in the first term of the sequence; this was shown to terminate at zero. It also implies that the reduction of the coefficient that multiplies the power 2, i.e. $(a_1 - 2) \cdot a_\chi^2$, has started. Hence, the sequence progresses up to a term in which such a coefficient is zero:

$$g_a(m, a_\delta) = (a_1 - 1) \cdot a_\delta^{a_1 - 1} + (a_1 - 1) \cdot a_\delta^{a_1 - 2} + \cdots + (a_1 - 1) \cdot a_\delta^4 + (a_1 - 1) \cdot a_\delta^3 - 1,$$

where $a_\delta = a_{\delta-1} + d_\delta$. Here, $(a_1 - 1) \cdot a_\delta^3 - 1 = (a_1 - 2) \cdot a_\delta^3 + a_\delta^3 - 1$. Thus, Lemma 3.4 will repeat the pattern, commencing at power 3. From the progression of the sequence up to the term $g_a(m, a_\delta)$, it is possible to deduce that, if a given power has been shown to terminate at zero, the same fate applies to that power multiplied by any coefficient (always smaller than the base). A further deduction is that, once the sequence reaches a given power in the absence of smaller ones (e.g. $a_\delta^3 - 1$ in the term $g_a(m, a_\delta)$ with all powers smaller than 3 having been eliminated), the application of Lemma 3.4 to that power implies that it will also be eliminated. Consequently, the pattern of elimination repeats until a term is reached such that

$$g_a(m, a_\varepsilon) = (a_1 - 1) \cdot a_\varepsilon^{a_1 - 1} - 1 = (a_1 - 2) \cdot a_\varepsilon^{a_1 - 1} + a_\varepsilon^{a_1 - 1} - 1,$$

where $a_\varepsilon = a_{\varepsilon-1} + d_\varepsilon$, and all powers smaller than $a_1 - 1$ have been eliminated. Lemma 3.4 will apply now to $a_\varepsilon^{a_1 - 1} - 1$, once more generating powers known to terminate at zero. The coefficient that multiplies this last power also begins to reduce in value. Therefore, it can be concluded that the complete sequence terminates at zero. □

Although the increases ($d$) in base value, as well as the starting base ($a_0$) of the sequence, take completely arbitrary values, as shown these do not impede the formulation of a satisfactory proof for the first part of the induction argument. On the contrary, this arbitrariness is essential to prove the second part of the argument, which deals with multiple levels of exponentiation.

*Step ii)* $m = a \Uparrow k_i$, and $m = a \Uparrow (k_i + 1)$. The second part of the induction process assumes that the generalized Goodstein sequence terminates at zero when $m = a \Uparrow k_i$. Thus, the requirement is to prove that the same can be stated for the next level of exponentiation, i.e. when $m = a \Uparrow (k_i + 1)$. Since $a \Uparrow (k_i + 1) = a^{(a \Uparrow k_i)}$, using the symbols $w_1$ and $v_1$ to signify $w_1 = (a_1 \Uparrow k_i) - 1$ and $v_1 = a_1 \Uparrow (k_i - 1) - 1$, respectively, and starting from $m = a_0 \Uparrow (k_i + 1)$, it can be written that

(5.1) $$g_a(m,a_1) = (a_1 \Uparrow (k_i+1)) - 1 = (a_1-1) \cdot a_1^{w_1} + (a_1-1) \cdot a_1^{w_1-1} + \cdots$$
$$\cdots + (a_1-1) \cdot a_1^{v_1} + \cdots + (a_1-1) \cdot a_1^{a_1} + \cdots + (a_1-1) \cdot a_1^{1} + (a_1-1) \cdot a_1^{0},$$

where $a_1 = a_0 + d_1$. The exponent of the first power in the summation is given by Lemma 3.4 as

$$w_1 = (a_1 \Uparrow k_i) - 1 = (a_1-1) \cdot a_1^{v_1} + (a_1-1) \cdot a_1^{v_1-1} + \cdots + (a_1-1) \cdot a_1^{1} + (a_1-1) \cdot a_1^{0}.$$

Significantly, the '$(k_i+1)$-level' of exponentiation, initially present in $a_1 \Uparrow (k_i+1)$, has already been eliminated in $g_a(m,a_1)$ since, due to Lemma 3.4, such a level can no longer appear in any of the powers in the summation that make the term (5.1). This reduction in the level of exponentiation plays a major role in the termination of the sequence[2].

There is a very important consideration to apply to (5.1), and repeatedly thereafter, that holds for any complete normal form composed of the sum of more than one power (and has effectively been used in the first part of the induction process); as a direct consequence of Lemma 4.1, if the largest power in the summation terminates at zero, all other (and smaller) powers in the form terminate at zero. Furthermore, the generalized Goodstein sequence for the whole summation terminates at zero. Because of the arbitrary value of the base, the initial assumption (i.e. that the largest power in the summation terminates at zero) maintains its validity while the progress of the sequence eliminates the smaller powers. The critical consequence of this analysis is that, in order to prove that a complete normal form terminates at zero, it is sufficient to prove it only for the largest power in the form. Thus, in the case of $m = a_0 \Uparrow (k_i+1)$, it will be sufficient to prove that the first power in (5.1) terminates at zero, i.e. $(a_1-1) \cdot a_1^{w_1}$. But also, given that $(a_1-1) \cdot a_1^{w_1} = (a_1-2) \cdot a_1^{w_1} + a_1^{w_1}$ (implying that the coefficient can be discounted), the requirement for proof of termination is reduced to just $a_1^{w_1}$.

In order to prove that the power $a_1^{w_1}$ terminates at zero, this power can be taken as the starting term of a new generalized Goodstein sequence such that

(5.2) $$a_2^{w_2} - 1 = (a_2-1) \cdot a_2^{w_2-1} + (a_2-1) \cdot a_2^{w_2-2} + \cdots + (a_2-1) \cdot a_2^{1} + (a_2-1) \cdot a_2^{0},$$

where $a_2 = a_1 + d_2$. In (5.2), $w_2$ is used to signify the complete normal form that results from replacing $a_1$ by $a_2$ in the representation of $w_1$. Applying the same analysis described for (5.1), it is possible to conclude that the termination of $a_1^{w_1}$ will be sufficiently proved by demonstrating only that $a_2^{w_2-1}$ terminates at zero.

As before, taking the latest power as the starting term of a new sequence, it can be written that

(5.3) $$a_3^{(w_2-1)_3} - 1 = (a_3-1) \cdot a_3^{(w_2-1)_3-1} + \cdots + (a_3-1) \cdot a_3^{1} + (a_2-1) \cdot a_3^{0},$$

where $a_3 = a_2 + d_3$. In (5.3), $(w_2-1)_3$ is used to signify the complete normal form that results from replacing $a_2$ by $a_3$ in the representation of $w_2-1$. Likewise, the termination of $a_2^{w_2-1}$ can be proved by demonstrating that $a_3^{(w_2-1)_3-1}$ terminates at zero.

---

[2] The exponents of all the powers in the term (5.1) are represented by complete normal forms to base $a_1$ that are smaller that $a_1 \Uparrow k_i$. Consequently, if generalized Goodstein sequences for the exponents themselves were to be considered, the assumption made would imply that all these sequences terminate at zero. Ultimately, this is the reason why the proof succeeds.

The above process can be continued repeatedly to generate a sequence of nested conditionals, where the proof of termination for one power depends on a similar proof for the subsequent power. The bases of the complete normal forms that represent the powers increase incrementally in value, as for any generalized Goodstein sequence. Thus, the sequential list of the exponents of the powers considered can be written as

$$(5.4) \qquad a_0 \Uparrow k_i \, , \, w_1 \, , \, w_2-1 \, , \, (w_2-1)_3-1 \, , \, ((w_2-1)_3-1)_4-1 \, , \, \ldots$$

The crucial observation is that the sequence (5.4) is the generalized Goodstein sequence for $m = a_0 \Uparrow k_i$. The assumption made for the induction step is that this sequence terminates at zero. Therefore, the sequence of nested conditionals constructed as described progresses to a power of zero exponent, which trivially terminates at zero. Since the last condition in the nested sequence is satisfied, the implication is that all the other conditions are equally satisfied, hence concluding that the initial power, $m = a_0 \Uparrow (k_i+1)$, terminates at zero, as required. □

*Step iii*) As indicated, the verification of *i*) and *ii*) implies that the sequence $g_a(m,a_n)$ terminates at zero when $m = a \Uparrow k$, $\forall a,k \in \mathbb{N}$. This conclusion, together with Lemma 3.3 and Lemma 4.1, implies also that the generalized Goodstein sequence $g_a(m,a_n)$ terminates at zero $\forall m \in \mathbb{N}$. □

Having proved Theorem 5.1, this result can be extended immediately to the more specific case of the original Goodstein's Theorem:

**Corollary 5.2 (Goodstein's Theorem).** *The Goodstein sequences of all numbers $m \in \mathbb{N}$ eventually terminate at zero:*

$$\forall m \in \mathbb{N}, \, \exists n \in \mathbb{N} \, (g(m,n) = 0).$$

*Proof.* The generalized Goodstein sequence defined by (4.1) is simplified to the Goodstein sequence defined by (2.3) by making $a_0 = 2$ and $d_n = 1$, $\forall n \in \mathbb{N}$. Hence, since (2.3) is only one example of (4.1), the proof of Theorem 5.1 doubles as proof for this corollary. □

## 6. Concluding remarks

The proof of Goodstein's Theorem reported here uses strictly finite means of deduction and operates exclusively with natural numbers. Thus, its construction lies entirely within the domain of first-order arithmetic (for example, PA). Since the generalized Goodstein sequence (4.1) and the super-exponential function (3.3) are both defined recursively, it should also be clear that the induction process used in the proof of Theorem 5.1 constitutes an instance of the first-order induction schema [10]. As explained in Section 1, Paris and Kirby established a firm connection between the inability to prove Goodstein's Theorem in PA and the axiomatic principles of ZFC [8,10]. Therefore, the proof described in this

12 J. A. PEREZarticle strongly challenges the integrity of ZFC. By providing an elementary argument that proves the truth of Goodstein's Theorem in first-order arithmetic, Theorem 5.1 and Corollary 5.2 establish a clear contradiction that implies the inconsistency of ZFC. The cause of this inconsistency has yet to be determined, but may lie in the way in which ZFC deals with mathematical infinities.

## References

[1] Gentzen G. *Die Widerspruchsfreiheit der reinen Zahlentheorie*. Math. Ann. 1936; 112: 493-565.
[2] Gentzen G. *Neue Fassung des Widerspruchsfreiheitsbeweises für die reine Zahlentheorie*. Forschungen zur Logik. 1938; 4: 19-44.
[3] Gödel, K. *Über formal unentscheidbare Sätze der* Principia Mathematica *und verwandter Systeme I*. Monatsh. Math. Physik. 1931; 38: 173-198.
[4] Goldrei D. *Classic Set Theory for Guided Independent Study*. Boca Raton: Chapman & Hall/CRC Press, 1998.
[5] Goodstein, R. L. *On the Restricted Ordinal Theorem*. J. Symb. Logic. 1944; 9: 33-41.
[6] Hawking S. *God Created the Integers: The Mathematical Breakthroughs that Changed History*. London: Penguin Books, 2006.
[7] Herrlich, H. *Lecture Notes in Mathematics: The Axiom of Choice*. Princeton: Princeton University Press, 2006.
[8] Kirby, L. and Paris, J. *Accessible independence results for Peano arithmetic*. Bull. London Math. Soc. 1982; 14: 285-293.
[9] Potter, M. *Set Theory and its Philosophy*. Oxford: Oxford University Press, 2004.
[10] Smith P. *An Introduction to Gödel's Theorems*. Cambridge: Cambridge University Press, 2007.
Juan A Perez. Berkshire, UK.
*E-mail address*: jap717@juanperezmaths.com